%

%
\input amstex
\documentstyle{amsppt}
\loadbold
\def\cstar{$C^*$-algebra}
\def\esg{$E_0$-semigroup}

\def\<{\left<}										
\def\>{\right>}

\magnification=\magstep 1

\topmatter
\title The index of a quantum dynamical semigroup
\endtitle
%


\author William Arveson
\endauthor

\affil Department of Mathematics\\
University of California\\Berkeley CA 94720, USA
\endaffil

\date 4 July 1996
\enddate
\thanks This research was supported by
NSF grant DMS95-00291
\endthanks
\keywords von Neumann algebras, automorphism groups, 
\esg s, minimal dilations, completely positive maps
\endkeywords
\subjclass
Primary 46L40; Secondary 81E05
\endsubjclass
\abstract 
A numerical index is introduced for semigroups 
of completely positive maps of $\Cal B(H)$ which 
generalizes the index of \esg s.   
It is shown that the index of a unital 
completely positive semigroup agrees with the
index of its dilation to an \esg , provided that 
the dilation is {\it minimal}.  
\endabstract

\endtopmatter
%

\document

\subheading{Introduction}
We introduce a numerical index for semigroups 
$P=\{P_t: t\geq 0\}$ of normal completely positive 
maps on the algebra $\Cal B(H)$ of all
bounded operators on a separable Hilbert space $H$.  
This index is defined 
in terms of basic structures associated with 
$P$, and generalizes the index of \esg s.  In the 
case where $P_t(\bold 1) = \bold 1$, $t\geq 0$, we 
show that the index of $P$ agrees with the index of 
its minimal dilation to an \esg.  

The key ingredients are the 
existence of the covariance function (Theorem 2.6),
the relation between units of $P$ and units of its
minimal dilation (Theorem 3.6), and the mapping 
of covariance functions (Corollary 4.8).  No 
examples are discussed here, but another paper 
is in preparation \cite{5}.

\subheading{1.  The metric operator space of a completely positive map}

We consider the real vector space of all normal 
linear maps $L$ of $\Cal B(H)$ into itself which are {\it symmetric\/} 
in the sense that $L(x^*) = L(x)^*$, $x\in \Cal B(H)$.  For two such maps
$L_1$, $L_2$ we write $L_1\leq L_2$ if the difference $L_2-L_1$ is 
completely positive.  Every operator $a\in \Cal B(H)$ gives rise to
an elementary completely positive map $\Omega_a$ by way of 
$$
\Omega_a(x) = axa^*, \qquad x\in \Cal B(H).  
$$

\proclaim{Definition 1.1}  For every completely positive map 
$P$ on $\Cal B(H)$ we write $\Cal E_P$ for the set of all operators
$a\in \Cal B(H)$ for which there is a positive constant $k$ 
such that
$$
\Omega_a \leq k P.  
$$
\endproclaim

\noindent 
In this section we collect some elementary observations which imply that
$\Cal E_P$ is a vector space inheriting a natural inner product
with respect to which it is a complex Hilbert space.  Thus, every 
normal completely positive map is associated with a Hilbert space
of operators which, as we will see, ``implement" the mapping.  
The properties of these Hilbert spaces of 
operators will be fundamental to our methods in the sequel.  

Because of Stinespring's theorem, every normal completely positive map
$P$ of $\Cal B(H)$ into itself can be represented in the form 
$$
P(x) = V^*\pi(x)V, \qquad x\in \Cal B(H),\tag{1.2}
$$
where $\pi$ is a representation of $\Cal B(H)$ on some Hilbert space
$H_\pi$ and $V:H\to H_\pi$ is a bounded operator.  We may always assume 
that the pair $(V,\pi)$ is {\it minimal} in the sense that $H_\pi$ 
is spanned by the set of vectors $\{\pi(x)V\xi: x\in \Cal B(H), \xi\in H\}$,   
and in that case we have $\pi(\bold 1) = \bold 1$ and $V^*V = P(\bold 1)$.   
Two minimal pairs $(V, \pi)$ and $(\tilde V, \tilde\pi)$ for $P$ 
are {\it equivalent} in the sense that there is a (necessarily unique) 
unitary operator $W: H_\pi \to H_{\tilde\pi}$ such that
$$
\align
WV &= \tilde V, \qquad \text{and} \tag{1.3.a}\\
W\pi(x) &= \tilde\pi(x)W, \qquad x\in \Cal B(H).  \tag{1.3.b}
\endalign
$$

Now since $P$ is normal, the representation $\pi$ occurring 
in any minimal pair $(V, \pi)$ is necessarily a normal representation
of $\Cal B(H)$ and is therefore unitarily equivalent to a representation
of the form 
$$
\pi(x) = x\oplus x\oplus  \dots,
$$
acting on a direct sum $H^n$ of $n$ copies of $H$, $n$ being a cardinal 
number which is countable because $H$ is separable.  Thus we may always assume
that a minimal pair $(V, \pi)$ consists of a representation of this form
and that $V: H\to H^n$ has the form 
$$
V\xi = (v_1^*\xi, v_2^*\xi, \dots)
$$
where $v_1, v_2, \dots$ is a sequence of bounded operators on $H$. 
Notice that the components of $V$ are the adjoints of the operators
$v_k$; this is essential in order for the operator multiplication 
to be properly related to the spaces $\Cal E_P$ associated with 
completely positive maps $P$ (see Theorem 1.12).   
After unravelling the formula (1.2) one finds that these operators 
satisfy
$$
P(x) = \sum_{n\geq 1} v_nxv_n^*
$$
the sum on the right converging weakly because of the condition
$$
\|v_1^*\xi\|^2 + \|v_2^*\xi\|^2 + \dots = \|V\xi\|^2 < \infty, \qquad \xi\in H. 
\tag{1.4} 
$$

Finally, the minimality condition on $(V,\pi)$ implies that the only
operator $c$ in the commutant of $\pi(\Cal B(H))$ satisfying 
$cV = 0$ is $c=0$.  Considering the matrix representation of operators 
in the commutant of $\pi(\Cal B(H))$, we find that this condition translates
into the somewhat more concrete ``linear independence" condition on 
the sequence $v_1, v_2, \dots$:
$$
(\lambda_1, \lambda_2, \dots) \in \ell^2,\quad \sum_k\lambda_k v_k = 0 \implies
\lambda_1 = \lambda_2 = \dots = 0.  \tag{1.5}
$$
Notice that the series in (1.5) is strongly convergent, since it 
represents the composition of $V^*:H^n\to H$ with the operator 
$\xi\in H \mapsto (\lambda_1\xi, \lambda_2\xi,\dots)\in H^n$.  
Conversely, if we start with an arbitrary sequence $v_1, v_2, \dots$ of 
operators in $\Cal B(H)$ for which (1.4) and the ``linear independence" 
condition (1.5) are satisfied, then 
$$
P(x) = \sum_k v_k x v_k^*
$$
defines a normal completely positive linear map on $\Cal B(H)$.  
If we define $V: H \to H^n$ and $\pi:\Cal B(H)\to \Cal B(H^n)$ by
$$
\align
V\xi &= (v_1^*\xi, v_2^*\xi, \dots), \tag{1.6.a}\\
\pi(x) &= x\oplus x\oplus \dots, \tag{1.6.b}
\endalign
$$ 
then $(V,\pi)$ is a minimal Stinespring pair $(V, \pi)$ for $P$.  

We now reformulate these observations in a coordinate-free
form which is more useful for our purposes below.  

\proclaim{Proposition 1.7}
Let $P(x) = V^*\pi(x)V$ be a minimal Stinespring representation 
for a normal completely positive map $P$ of $\Cal B(H)$, and let
$$
\Cal S = \{T\in \Cal B(H, H_\pi) : Tx = \pi(x) T, x\in \Cal B(H)\}
$$
be the intertwining space for $\pi$ and the identity representation.  
For any two operators $T_1, T_2 \in \Cal S$, $T_2^*T_1$ is a scalar 
multiple of the identity of $\Cal B(H)$, and 
$$
\<T_1, T_2\>\bold 1 = T_2^*T_1
$$
defines an inner product on $\Cal S$ with respect to which it is 
a Hilbert space in which the operator norm coincides with the Hilbert 
space norm.  

The linear mapping $T\in \Cal S \to V^*T\in \Cal B(H)$ is injective 
and has range $\Cal E_P$.  $\Cal E_P$ is  
a Hilbert space with respect to the inner product defined 
by pushing forward the inner product of $\Cal S$.  
\endproclaim
\demo{proof}
We merely sketch the argument, which is part of the folklore of 
representation theory.  The first paragraph is completely straightforward.  
For example, if $T_1, T_2\in \Cal S$ then $T_2^*T_1$ must be a scalar 
multiple of the identity on $H$ because for every $x\in \Cal B(H)$ we have
$$
T_2^*T_1 x = T_2^*\pi(x)T_1 = x T_2^*T_1.  
$$

Now let $a$ be an operator of the form
$a = V^*T$, $T\in \Cal S$.  We claim that $a$ belongs to $\Cal E_P$.  
Indeed, for every $x\in \Cal B(H)$ we have 
$$
\Omega_a(x) = axa^* = V^*TxT^*V = V^*\pi(x)TT^*V.  
$$
Since $TT^*$ is a bounded positive operator in the commutant of 
$\pi(\Cal B(H))$, the operator $C = (\|T\|^2\bold 1 - TT^*)^{1/2}$
is positive, commutes with $\pi(\Cal B(H))$, and 
the preceding formula implies that the operator mapping  
$$
x\in \Cal B(H) \mapsto \|T\|^2 P(x) - \Omega_a(x) = V^*C\pi(x)CV  
$$
is completely positive.  Hence $a\in \Cal E_P$.  

The map $T\to V^*T$ is injective because it is linear, and because if an
operator $T\in \Cal S$ satisfies $V^*T = 0$ then for every $x\in \Cal B(H)$ 
and every $\xi\in H$ 
$$
T^*\pi(x)V\xi = xT^*V\xi = x(V^*T)^*\xi = 0.  
$$ 
Hence $T^*=0$ because $H_\pi$ is spanned by $\pi(\Cal B(H))H$, 
hence $T=0$.  

Finally, let $a$ be an arbitrary element in $\Cal E_P$ and choose
a positive constant $k$ such that $\Omega_a \leq k P$.  We may find 
an operator $T\in \Cal S$ which maps to $a$ as follows.  For any $n\geq 1$, 
any operators $x_1, x_2, \dots x_n\in \Cal B(H)$ 
and any vectors $\xi_1, \xi_2, \dots \xi_n\in H$ we have
$$
\align
\|\sum_{k=1}^n x_ka^*\xi_k\|^2 &= \sum_{j,k = 1}^n\<\Omega_a(x_k^*x_j)\xi_j,\xi_k\>
\leq k \sum_{k, j=1}^n\<V^*\pi(x_k^*x_j)V\xi_j, \xi_k\> \\
&=k\|\sum_{k=1}^n\pi(x_k)V\xi_k\|^2.  
\endalign
$$
\noindent
Thus there is a unique bounded operator $L: H_\pi = [\pi(\Cal B(H)VH]\to H$
which satisfies $L(\pi(x)V\xi) = xa^*\xi$ for every $x\in \Cal B(H)$,
$\xi\in H$.  Taking $T=L^*$ we find that $T\in \Cal S$ and $a^*=T^*V$,
hence $a = V^*T$ \qed
\enddemo

\remark{Remark 1.8}
To reiterate, the inner product in $\Cal E_P$ is defined as follows.  
Pick $a, b\in \Cal E_P$.  
Then there are unique operators $S, T\in \Cal S$ such 
that $a = V^*S$, $b=V^*T$, and $\<a,b\>$ is defined by
$$
\<a,b\>\bold 1 = T^*S.  
$$

In more concrete terms, choose a minimal Stinespring representation
$P(x) = V^*\pi(x)V$ where $\pi$ is a representation on $H^n$ and 
$V: H\to H^n$ is of the form $V\xi = (v_1\xi, v_2\xi, \dots)$, the
sequence of operators $v_1, v_2, \dots\in \Cal B(H)$ satisfying 
conditions (1.4) and (1.5).  Then $\{v_1^*, v_2^*, \dots\}$ is 
an orthonormal basis for the Hilbert space structure of 
$\Cal E_P$ and thus $\Cal E_P$ consists precisely of
all operators $a$ of the form
$$
a = a(\lambda) = \lambda_1 v_1^* + \lambda_2 v_2^* +\dots,
$$
where $\lambda = (\lambda_1, \lambda_2, \dots)$ is an arbitrary sequence 
in $\ell ^2$.  The sequence $\lambda\in\ell ^2$ is uniquely determined
by the operator $a(\lambda)$, and the inner product in $\Cal E_P$ satisfies
$$
\<a(\lambda), a(\mu)\> = \sum_{k\geq 1}\lambda_k\bar\mu_k.  
$$
\endremark

\proclaim{Definition 1.9}
A {\bf metric operator space} is a pair $(\Cal E,\<\cdot, \cdot\>)$ 
consisting of a complex 
linear subspace $\Cal E$ of $\Cal B(H)$ together with an 
inner product $u, v\in \Cal E\mapsto \<u,v\>\in \Bbb C$ with respect 
to which $\Cal E$ is a separable Hilbert space 
which has the following property: 
if $e_1, e_2, \dots$ is an orthonormal basis for $\Cal E$ then for 
any $\xi \in H$ we have
$$
\|e_1^*\xi\|^2 + \|e_2^*\xi\|^2 + \dots < \infty.  \tag{1.10}
$$
\endproclaim

\remark{Remarks}
The above discussion shows how, starting with a normal completely 
positive map $P$ of $\Cal B(H)$ into itself, we associate with
$P$ in an invariant way a metric operator space $\Cal E_P$.  
This metric operator space has the property that if we pick 
an arbitrary orthonormal basis $e_1, e_2, \dots$ for $\Cal E$ 
then we recover the map $P$ as follows,
$$
P(x) = \sum_k e_k x e_k^*, \qquad x\in \Cal B(H), \tag{1.11}
$$
the sum on the right being independent 
of the particular choice of basis.  Conversely, starting with an
arbitrary metric operator space $\Cal E$ we may define a 
unique completely positive map $P$ by the formula (1.11), and
thus {\it we have a bijective correspondence $P\leftrightarrow \Cal E$
between normal completely positive maps and metric operator 
spaces}.  

Metric operator spaces offer several advantages over the 
Stinespring representation in describing normal completely positive 
maps on $\Cal B(H)$, and it is appropriate to briefly discuss these 
issues here.  For example, suppose we start with such a map $P$ 
with metric operator space $\Cal E$.  We may use the 
inner product on $\Cal E$ to define an inner product on the 
tensor product of vector spaces $\Cal E\odot H$, and after 
completion we obtain a Hilbert space $\Cal E\otimes H$.  
The natural multiplication map $M: \Cal E\odot H\to H$ 
defined by $M(v\otimes\xi)=v\xi$ extends uniquely to a 
bounded operator from $\Cal E\otimes H$ to $H$, which we 
denote by the same letter $M$.  To see that, 
choose an orthonormal basis $e_1,e_2,\dots$ for 
$\Cal E$ and define a (necessarily bounded) operator 
$V:H\to\Cal E\otimes H$ by 
$$
V\xi=\sum_ke_k\otimes e_k^*\xi.  
$$
A direct computation then shows that 
$$
\<M(v\otimes\xi),\eta\>_H = 
\<v\otimes\xi,V\eta\>_{\Cal E\otimes H},
$$ 
and hence $M=V^*$.  In particular, the operator $V$ is 
independent of the particular choice of basis, 
and represents ``comultiplication".  
Moreover, if we define a normal 
representation $\pi: \Cal B(H)\to\Cal B(\Cal E\otimes H)$
by 
$$
\pi(x)=\bold 1_{\Cal E}\otimes x
$$
then one finds that $(V,\pi)$ is a {\it minimal} 
Stinespring representation for $P$.  We conclude that 
with every normal completely positive map $P$ there is a
natural way of picking out a concrete minimal Stinespring 
pair $(V,\pi)$ for $P$: one computes the metric operator 
space $\Cal E$ associated with $P$, takes 
$V:H\to\Cal E\otimes H$ to be comultiplication and
takes $\pi$ as above.  

More significantly, notice that the Stinespring 
representation of normal completely positive 
maps does not behave well with respect 
to composition.  For example, if we have two such maps 
$P_k: \Cal B(H)\to\Cal B(H)$, and we consider their 
respective minimal Stinespring pairs $(V_k, \pi_k)$, then
there is no natural way to combine $(V_1,\pi_1)$ with 
$(V_2, \pi_2)$ to obtain a Stinespring pair for the 
composition $P_1P_2$, much less 
a minimal one.  The description 
of such maps in terms of metric operator spaces is 
designed to deal efficiently with compositions.  
The following result implies that 
the metric operator space of $P_1P_2$
is spanned (as a Hilbert space) by the set of all 
operator products $\Cal E_1\Cal E_2$, 
$\Cal E_k$ denoting the space associated with $P_k$.  As 
we will see in the sequel, this is a critical feature 
when dealing with semigroups.  
\endremark
\proclaim{Theorem 1.12}
Let $\Cal E_1$ and $\Cal E_2$ be metric operator spaces 
with corresponding completely positive maps $P_k = P_{\Cal E_k}$, 
and let $P_1P_2$ denote the composition.  
Let $\Cal E_1\otimes \Cal E_2$ be the tensor product of
Hilbert spaces.  Then $\Cal E_{P_1P_2}$ contains the 
set of  all operator products 
$\{uv: u\in \Cal E_1, v\in \Cal E_2\}$ and 
there is a unique bounded linear operator
$M:\Cal E_1\otimes\Cal E_2\to\Cal E_{P_1P_2}$ satisfying
$$
M(u\otimes v) = uv.  \qquad u\in \Cal E_1, v\in \Cal E_2, \tag{1.13}
$$
The adjoint of $M$ is an isometry
$$
M^*: \Cal E_{P_1P_2}\hookrightarrow \Cal E_1\otimes\Cal E_2 
$$
whose range is a (perhaps proper) closed subspace of 
$\Cal E_1\otimes\Cal E_2$.  
\endproclaim

\remark{Remarks}
We refer to the adjoint $M^*$ of the operator $M$
defined by (1.13) as {\it comultiplication}.  
Since comultiplication is an 
isometry, it follows that the range of the multiplication
operator $M$ is all of $\Cal E_{P_1P_2}$, and hence 
$\Cal E_{P_1P_2}$ is spanned by the set 
of products $\Cal E_1\Cal E_2$.  

Theorem 1.12 asserts that comultiplication gives rise to 
a natural identification of $\Cal E_{P_1P_2}$ with a 
closed subspace of $\Cal E_1\otimes\Cal E_2$.  Equivalently, 
the polar decomposition of the multiplication operator $M$ 
has the form $M = UQ$, where 
$Q\in \Cal B(\Cal E_1\otimes\Cal E_2)$ is
the projection onto this subspace and 
$U\in \Cal B(\Cal E_1\otimes\Cal E_2,\Cal E_{P_1P_2})$ is 
a partial isometry with $U^*U = Q$, 
$UU^* = \bold 1_{\Cal E_{P_1P_2}}$.  
\endremark

\demo{proof of Theorem 1.12}
We find a Stinespring representation of $P_1P_2$ in terms 
of $\Cal E_1$ and $\Cal E_2$ as follows.  Consider the 
Hilbert space $K = \Cal E_1\otimes\Cal E_2\otimes H$, and 
the representation $\pi$ of $\Cal B(H)$ on $K$ defined by
$$
\pi(x) = \bold 1_{\Cal E_1}\otimes\bold 1_{\Cal E_2}\otimes x.  
$$
Choose an orthonormal basis $u_1, u_2,\dots$ for $\Cal E_1$ 
(resp. $v_1, v_2, \dots$ for $\Cal E_2$) and define an 
operator $V: H\to K$ by 
$$
V\xi = \sum_{i,j} u_i\otimes v_j\otimes v_j^*u_i^*\xi, 
\qquad \xi\in H.  
$$
It is clear that $V$ is bounded, since
$$
\|V\xi\|^2 = \sum_{i,j}\|v_j^*u_i^*\xi\|^2 = 
\sum_{i,j}\<u_iv_jv_j^*u_i^*\xi,\xi\> = 
\<P_1(P_2(\bold 1))\xi,\xi\>,
$$
and in fact $V^*V = P_1P_2(\bold 1)$.  A similar calculation 
shows that 
$$
V^*\pi(x) V = P_1P_2(x), \qquad x\in \Cal B(H).  
$$

However, $(V,\pi)$ is not necessarily a {\it minimal} 
Stinespring pair.  In order to arrange minimality, 
consider the subspace $K_0\subseteq K$ defined by
$$
K_0 = [\pi(x)V\xi: x\in \Cal B(H), \xi\in H].  
$$
Since $K_0$ is invariant under the range of $\pi$ its production
belongs to the commutant
$$
\pi(\Cal B(H))^\prime = 
\Cal B(\Cal E_1\otimes\Cal E_2)\otimes \bold 1_H,
$$
and hence there is a unique projection 
$Q\in \Cal B(\Cal E_1\otimes\Cal E_2)$ such that 
$$
P_{K_0} = Q\otimes \bold 1_H.  
$$
The corresponding subrepresentation 
$\pi_0$ obtained by restricting $\pi$ to $K_0$ gives rise
to a minimal Stinespring pair $(V,\pi_0)$ for $P_1P_2$.  

In order to calculate the metric operator space 
$\Cal E_{P_1P_2}$ we use Proposition 1.7 as follows.  Notice
that for every $\zeta\in\Cal E_1\otimes\Cal E_2$ we can 
define a bounded operator $X_\zeta: H\to K$ by
$$
X_\zeta \xi = \zeta\otimes \xi, \qquad \xi\in H.  
$$
It is clear that $X_\zeta a = \pi(a)X_\zeta$ for every 
$a\in \Cal B(H)$, and moreover every bounded operator 
$X\in \Cal B(H,K)$ satisfying $Xa = \pi(a)X$, $a\in \Cal B(H)$, 
has the form $X=X_\zeta$ for a unique 
$\zeta\in \Cal E_1\otimes\Cal E_2$.  The range of 
$X_\zeta$ is contained in $K_0=(Q\otimes\bold 1_H)K$ if and 
only if $\zeta$ belongs to the range of $Q$.  
Thus the intertwining space 
$$
\{X\in\Cal B(H,K_0):Xa = \pi_0(a)X, a\in \Cal B(H)\}
$$ 
for $\pi_0$ is 
$\{X_\zeta: \zeta\in Q(\Cal E_1\otimes\Cal E_2)\}$.

Now by Proposition 1.7, we have 
$$
\Cal E_{P_1P_2} = 
\{V^*X_\zeta: \zeta\in Q(\Cal E_1\otimes\Cal E_2)\}, 
$$
and the inner product of two operators $T_k = V^*X_{\zeta_k}$, 
$k=1,2$ in $\Cal E_{P_1P_2}$ is given by 
$$
\<T_1,T_2\>_{\Cal E_{P_1P_2}}\bold 1_H = X_{\zeta_2}^*X_{\zeta_1} 
= \<\zeta_1,\zeta_2\>\bold 1_H,
$$
$\zeta_k\in Q(\Cal E_1\otimes\Cal E_2)$.  

Accordingly, we have defined a unitary 
operator $U: Q(\Cal E_1\otimes\Cal E_2)\to \Cal E_{P_1P_2}$ 
by 
$$
U\zeta = V^*X_\zeta.  
$$
It remains to show that the bounded operator 
$M:\Cal E_1\otimes\Cal E_2\to \Cal E_{P_1P_2}$ defined by
$$
M = UQ
$$
represents multiplication in the sense that 
$M(u\otimes v) = uv$ for any 
$u\in \Cal E_1$ and $v\in \Cal E_2$.  To see that, write 
$$
M(u\otimes v) = UQ(u\otimes v) = V^*X_{Q(u\otimes v)} =
V^*(Q\otimes \bold 1_H)X_{u\otimes v} = V^*X_{u\otimes v}, 
$$
the last equality following from the fact that 
$Q\otimes \bold 1_H V = P_{K_0}V = V$.  Thus for 
$\xi$, $\eta\in H$ we have 
$$
\align
\<M(u\otimes v)\xi,\eta\> &= \<V^*X_{u\otimes v}\xi,\eta\> 
=\sum_{i,j}\<u\otimes v\otimes\xi,
u_i\otimes v_j\otimes v_j^*u_i^*\eta\> \\
&=\sum_{i,j}\<u,u_i\>\<v,v_j\>\<\xi,v_j^*u_i^*\eta\> =
\sum_{i,j}\<u,u_i\>\<v,v_j\>\<u_iv_j\xi,\eta\>. 
\endalign 
$$
The term on the right is $\<uv\xi,\eta\>$ because 
$\sum_i\<u,u_i\>u_i = u$ and $\sum_j\<v,v_j\>v_j=v$\qed 
\enddemo

\remark{Remark 1.14}
Finally, we call attention to the special 
case in which $P$ is a normal $*$-endomorphism, that
is, a normal completely positive map 
for which $P(xy) = P(x)P(y)$ for all $x, y$.  
We do not assume that $P(\bold 1) = \bold 1$, 
but of course $P(\bold 1)$
must be a self-adjoint projection.  
In this case a minimal Stinespring
representation $P = V^*\pi V$ is given 
by the pair $(V, \pi)$, where V is 
the orthogonal projection of $H$ onto 
$H_0 = P(\bold 1)H$ and $\pi(x)$ is 
the restriction of $P(x)$ to the 
invariant subspace $H_0$.  In this 
case a straightforward computation 
shows that $\Cal E_P$ reduces to 
the intertwining space 
$$
\Cal E_P = \{T\in \Cal B(H): P(x)T = T x, x\in \Cal B(H)\}, 
$$
and that the inner product on $\Cal E_P$ is defined by 
$$
\<T_1, T_2\> \bold 1 = T_2^* T_1, \qquad T_1, T_2 \in \Cal E_P.  
$$
\endremark

\subheading{2.  Numerical index}
Let $H$ be a separable Hilbert space and let 
$P = \{P_t: t\geq 0\}$ be a semigroup of normal completely positive 
maps of $\Cal B(H)$ into itself which is continuous in the sense 
that for every $x\in \Cal B(H)$ and every pair of vectors 
$\xi, \eta\in \Cal B(H)$, the function 
$t\in [0,\infty)\mapsto \<P_t(x)\xi,\eta\>$ is continuous.  We 
do not assume that $P_t$ preserves the unit, nor even that 
$\|P_t\|\leq 1$, but we do require that $P_0$ be the identity
map; equivalently, 
$$
\lim_{t\to 0}\<P_t(x)\xi,\eta\> = \<x\xi,\eta\>, 
\qquad x\in \Cal B(H), \xi,\eta\in H.  
$$
We refer 
to such a semigroup as a {\it CP semigroup}.  A CP semigroup $P$ is 
called {\it unital} if $P_t(\bold 1) = \bold 1$ for every $t\geq 0$, 
and {\it contractive} if $\|P_t\| = \|P_t(\bold 1)\| \leq 1$ for 
every $t\geq 0$.  

In this section we introduce a 
numerical index for arbitrary CP semigroups 
which generalizes the definition
of index of  \esg s \cite{1}.  While the definition 
and Theorem 2.6 below are very general, the reader should 
keep in mind that we are primarily 
interested in the case of {\it unital} CP semigroups.  

\proclaim{Definition 2.1}
Let $P$ be a CP semigroup acting on $\Cal B(H)$.  
A {\bf unit} of $P$ is a semigroup
$T = \{T_t: t\geq 0\}$ of bounded operators on $H$ which is strongly
continuous in the sense that
$$
\lim_{t\to 0}\|T_t\xi - \xi\| = 0, \qquad \xi\in H
$$
and for which there is a real constant $k$ such 
that for every $t > 0$, the operator mapping 
$\Omega_t(x) = T_txT_t^*$ satisfies 
$$
\Omega_t \leq e^{kt}P_t.  
$$
\endproclaim

\remark{Remark 2.2}
We write $\Cal U_P$ for the set of all units of $P$, and it will
be convenient to denote the metric operator spaces $\Cal E_{P_t}$ 
associated with the individual 
completely positive maps $P_t$ with the notation
$\Cal E_P(t)$, $t\geq 0$.  Notice that an operator 
semigroup $T = \{T_t: t\geq 0\}$ belongs to $\Cal U_P$ if and only
if a) $T_t \in \Cal E_P(t)$ for every $t > 0$ and b) the
Hilbert space norms $\<T_t, T_t\>$ of these elements of $\Cal E_P(t)$
satisfy the growth condition
$$
\<T_t, T_t\> \leq e^{kt}, \qquad t > 0.  \tag{2.3}
$$
\endremark

Of course, every \esg\ qualifies as a CP semigroup, and in this 
case remark (1.14) implies that Definition 2.1 agrees with the 
definition of unit for an \esg\ given in \cite{1}.  The only issue
here is the growth condition (2.3), which is not part of the definition
of unit for an \esg.  However, if $T = \{T_t: t>0\}$ is a unit for 
an \esg\ $P = \{P_t : t\geq 0\}$ then we have 
$$
\<T_t, T_t\> = e^{tc(T, T)}
$$
where $c: \Cal U_P\times \Cal U_P\to \Bbb C$ is the covariance 
function defined in \cite{1}, and hence the growth condition 
(2.3) is {\it automatic} for \esg s.  

Since there exist \esg s with no units whatsoever \cite{9} we 
must allow for the possibility that a CP semigroup may have 
no units.  However, assuming that $P$ is a CP semigroup 
for which $\Cal U_P \neq \emptyset$, we define a numerical 
index $d_*(P)$ in the following way.  Choose $S, T\in \Cal U_P$.  
Then for every $t>0$ the operators $S_t, T_t$ belong to the Hilbert
space $\Cal E_P(t)$ and we may consider their inner product
$$
\<S_t, T_t\>\in \Bbb C. \tag{2.4}
$$  
Notice that while the inner products (2.4) are computed in different
Hilbert spaces $\Cal E_P(t)$, there is no ambiguity in this notation 
so long as the variable $t$ is displayed.  We remark too 
that while neither semigroup $S$ nor $T$ can be the zero semigroup, 
it can certainly happen that $T_t = 0$ for certain positive values 
of $t$, and once $T_t$ is zero for some particular
value of $t$ then it is zero for all larger $t$ as well.  However, 
strong continuity at $t=0$ implies that for sufficiently small $t$,
both operators $S_t$ and $T_t$ are nonzero.  But even in this 
case, there is no obvious guarantee that the inner product 
$\<S_t, T_t\>$ is nonzero.  

Now fix $t>0$ and choose $S, T\in \Cal E_P(t)$.  For each finite 
partition 
$$
\Cal P = \{0 = t_0<t_1<\dots <t_n = t\}
$$ 
of the interval $[0,t]$ we define 
$$
f_\Cal P(S, T; t) = \prod_{k=1}^n\<S_{t_k-t_{k-1}}, T_{t_k-t_{k-1}}\>.  
\tag{2.5}
$$
If we consider the set of partitions of $[0, t]$ as an increasing 
directed set in the usual way then (2.5) defines a net of complex
numbers.  The definition of index depends on the following 
result, which will be proved later in this section.  

\proclaim{Theorem 2.6}Let $P = \{P_t: t\geq 0\}$ be a 
CP semigroup acting on $\Cal B(H)$, let $S$ and $T$ be units
of $P$, and define $f_{\Cal P}(S,T;t)$ as in (2.5).  
Then there is a (necessarily unique) complex number $c$ such that 
$$
\lim_\Cal P f_\Cal P(S, T; t) = e^{ct} 
$$
for every $t>0$.  
\endproclaim

We postpone the proof of Theorem 2.6 in order to discuss 
its immediate consequences.  
We will write $c_P(S, T)$ for the constant $c$ of 
Theorem 2.6.  Thus we have defined a bivariate function 
$$
c_P: \Cal U_P \times \Cal U_P \to \Bbb C, 
$$
which will be called the {\bf covariance function} of the CP 
semigroup $P$.  

\proclaim{Proposition 2.7}
The covariance function is conditionally positive definite in the sense
that if $T_1, T_2, \dots, T_n\in \Cal U_P$ and $\lambda_1, \lambda_1, 
\dots, \lambda_n$ are complex numbers satisfying 
$\lambda_1 + \lambda_2 + \dots + \lambda_n = 0$, then 
$$
\sum_{j,k = 1}^n \lambda_j\bar\lambda_k c_P(T_j, T_k) \geq 0.  
$$
\endproclaim
\demo{proof}
It suffices to show that for every fixed $t>0$, the function
$$
S, T \mapsto e^{tc_P(S, T)} 
$$
is positive definite \cite{8}.  Now for every positive
$\lambda$, $S, T\mapsto \<S_\lambda, T_\lambda\>$ 
is obviously a positive definite function.  Since a finite 
pointwise product of positive definite functions 
is a positive definite function, it follows that for each 
partition $\Cal P$ of $[0, t]$ the function 
$$
S, T \mapsto f_\Cal P(S, T; t)
$$
of (2.5) is positive definite.  Finally, since the limit of a
pointwise convergent net of positive definite functions is 
positive definite, we conclude from Theorem 2.6 that the
function 
$$
e^{tc_P(S,T)} = \lim_\Cal P f_\Cal P(S,T;t)
$$ 
must be a positive definite of $S$ and $T$ \qed
\enddemo

Now suppose that $P=\{P_t: t\geq 0\}$ is a CP semigroup
for which $\Cal U_P\neq \emptyset$.  
We may construct a Hilbert space $H_P$ out of
the conditionally postive definite function 
$c_P: \Cal U_P\times \Cal U_P\to \Bbb C$
in the same way as for \esg s.  More explicitly, 
on the vector space $V$ consisiting of all finitely 
nonzero functions $f:\Cal U_P\to\Bbb C$ satisfying  
$$
\sum_{T\in\Cal U_P}f(T)=0,
$$
one defines a positive semidefinite sesquilinear 
form 
$$
\<f,g\>=\sum_{S,T\in\Cal U_P}f(S)\overline{g(T)}c_P(S,T),   
$$
and the Hilbert space $H_P$ is obtained by completing 
the inner product space $V/N$, where $N$ is the subspace
$$
N=\{f\in V: \<f,f\>=0\}.  
$$
We define the 
{\bf index} of $P$ as the dimension of this Hilbert space
$$
d_*(P) = \dim(H_P).  
$$
For the principal class of examples in which 
$P$ is a {\it unital} CP semigroup Corollary 4.8 below
together with \cite{1, Proposition 5.2} implies that $H_P$ must 
be separable, so that $d_*(P)$ must take one of the values
$0, 1, 2, \dots, \aleph_0$.  

The exceptional case in which $\Cal U_P = \emptyset$ is handled in 
the same way as for \esg s; in that event we define 
$$
d_*(P) = 2^{\aleph_0}
$$
to be the cardinality of the continuum.  This convention of 
choosing an uncountable value for the index in the exceptional
case where there are no units allows for
the unrestricted validity of the addition formula
for tensor products
$$
d_*(P\otimes Q) = d_*(P) + d_*(P)  
$$
in the same way it does for \esg s.  

\demo{proof of Theorem 2.6}
Let $T_k = \{T_k(t): t\geq 0\}$, $k=1,2$, be units  of 
a fixed CP semigroup $P$.  Because 
each unit $T$ must satisfy a growth condition of the form 
$\<T(t),T(t)\>\leq e^{ct}$, $t>0$, we may rescale $T_1$ and 
$T_2$ with a factor of the form $e^{-c_kt}$ to achieve 
$$
\<T_k(t),T_k(t)\>\leq 1, \qquad t>0.  \tag{2.8}
$$
Notice that this rescaling does not affect either the existence 
of the limit of Theorem 2.6 or the exponential nature of its value, 
so it suffices to prove 2.6 in the presence of the normalization
(2.8).

For each partition 
$\Cal P = \{0=t_0<t_1<\dots<t_n=t\}$ of the interval $[0,t]$ 
we consider the $2\times 2$ matrix $A_\Cal P(t)$ whose $ij$th 
term is given by
$$
f_\Cal P(T_i,T_j;t) = \prod_{k=1}^n\<T_i(t_k-t_{k-1}),T_j(t_k-t_{k-1})\>. 
\tag{2.9} 
$$
(2.8) implies that $|f_\Cal P(T_i,T_j,;t)|\leq 1$; thus
we have a uniform bound
$$
\|A_\Cal P(t)\|\leq \text{trace}(A_\Cal P(t)^*A_\Cal P(t))^{1/2} \leq 2.  
$$

As in the proof of Proposition (2.7), each function 
$f_\Cal P(\cdot,\cdot; t)$ is positive definite; hence
$A_\Cal P(t)$ is a positive matrix.  We claim that in fact, 
$$
\Cal P_1\subseteq \Cal P_2\implies 
A_{\Cal P_1}(t)\leq A_{\Cal P_2}(t).
\tag{2.10}
$$
To see that, it is enough to consider the case where 
$\Cal P_2$ is obtained by adjoining a single point $\lambda$ to 
$\Cal P_1 = \{0=t_0<t_1<\dots<t_n=t\}$.  
Suppose that $t_{k-1}<\lambda<t_k$ for $k$ between
$1$ and $n$.  Note that $f_{\Cal P_2}(T_i,T_j;t)$ is obtained 
from $f_{\Cal P_1}(T_i,T_j;t)$ by replacing the $k$th term 
$\alpha_{ij} = \<T_i(t_k-t_{k-1}),T_j(t_k-t_{k-1})\>$ in the product (2.9) 
with the term 
$$
\beta_{ij} = \<T_i(\lambda-t_{k-1}),T_j(\lambda-t_{k-1})\>
\<T_i(t_k-\lambda),T_j(t_k-\lambda)\>.  
$$
Thus, the $ij$th term of $A_{\Cal P_2}(t)-A_{\Cal P_1}(t)$ 
has the form $(\beta_{ij}-\alpha_{ij})\gamma_{ij}$, where the 
$2\times 2$ matrix $(\gamma_{ij})$ is positive. 
Since the 
Schur product of two positive matrices is positive, it 
suffices to show that $(\beta_{ij}-\alpha_{ij})$ is a positive
$2\times 2$ matrix.  Now for any two complex numbers 
$\lambda_1,\lambda_2$ we have 
$$
\align
&\sum_{i,j=1}^2\lambda_i\bar\lambda_j\beta_{ij}-
\sum_{i,j=1}^2\lambda_i\bar\lambda_j\alpha_{ij} =\\
&\sum_{i,j=1}^2\lambda_i\bar\lambda_j
\<T_i(\lambda-t_{k-1}),T_j(\lambda-t_{k-1})\> 
\<T_i(t_k-\lambda),T_j(t_k-\lambda)\>-\\
&\sum_{i,j=1}^2\lambda_i\bar\lambda_j
\<T_i(t_k-t_{k-1}),T_j(t_k-t_{k-1})\> = \\
&\|\sum_i\lambda_i T_j(\lambda-t_{k-1}))\otimes T_i(t_k-\lambda)\|^2 -
\|\sum_i\lambda_iT_i(t_k-t_{k-1})\|^2.  
\endalign
$$
Because of the semigroup property we have 
$T_i(t_k-t_{k-1})=T_i(\lambda-t_{k-1})T_i(t_k-\lambda)$.  
Thus the last term of the preceding 
formula is nonnegative because of Theorem 1.12, 
which implies that multiplication 
$$
M: \Cal E_P(\lambda-t_{k-1})\otimes\Cal E_P(t_k-\lambda) \to 
\Cal E_P(t_k-t_{k-1})
$$
is a contraction.  This establishes (2.10).  

Since for fixed $t>0$, $\Cal P\mapsto A_\Cal P(t)$ is a 
uniformly bounded 
increasing net of positive operators, conventional wisdom implies 
that there is a unique positive operator $B(t)\in M_2(\Bbb C)$ 
such that 
$$
B(t) = \lim_\Cal P A_\Cal P(t).  
$$
Letting $b_{ij}(t)$ be the $ij$th entry of $B(t)$ we 
have the required limit (2.6),
$$
b_{ij}(t) = \lim_\Cal P f_{\Cal P}(T_i,T_j;t).  \tag{2.11}
$$

It remains to show that 
the functions $b_{ij}$ have the form 
$$
b_{ij}(t) = e^{tc_{ij}} \qquad t>0\tag{2.12}
$$
for some $2\times 2$ matrix $(c_{ij})$.  Now every pair 
$\Cal P$, $\Cal Q$ consisting of finite partitions of 
$[0,s]$ and $[0,t]$ respectively gives rise to a partition 
of $[0,s+t]$, simply by first listing the elements of $\Cal P$
and then listing the elements of $s + \Cal Q$.  This construction 
gives all partitions of $[0,s+t]$ which contain the 
point $s$.  Since the latter is a cofinal subset of 
all finite partitions
of $[0,s+t]$ it follows from (2.11) that we have 
$$
b_{ij}(s+t) = b_{ij}(s)b_{ij}(t), \qquad s,t>0.  
$$
Thus to prove (2.12) it is enough to show that the functions
$b_{ij}$ extend continuously to the origin in the following
sense
$$
\lim_{t\to 0+}b_{ij}(t) = 1.  
$$
The latter is an immediate consequence of the following 
two results.

\proclaim{Lemma 2.14} For $i,j=1$ or $2$ and $t>0$ we have 
$$
|b_{ij}(t) - \<T_i(t),T_j(t)\> |^2\leq 
(1-\<T_i(t),T_i(t)\>)(1-\<T_j(t),T_j(t)\>).  
$$
\endproclaim

\proclaim{Lemma 2.15}For $i,j = 1$ or $2$ we have 
$$
\lim_{t\to 0+}\<T_i(t),T_j(t)\> = 1.  
$$
\endproclaim

\demo{proof of Lemma 2.14}Fix $t>0$.  
Because of (2.11), it suffices to 
show that for every $i$ and $j$ and every finite 
partition $\Cal P = \{0=t_0<t_1<\dots<t_n=t\}$ 
of the interval $[0,t]$, we have 
$$
|f_\Cal P(T_i,T_j;t) - \<T_i(t),T_j(t)\>|^2 \leq 
(1-\<T_i(t),T_i(t)\>)(1-\<T_j(t),T_j(t)\>).  \tag{2.16}
$$
Consider the vectors 
$u_i\in \Cal E_P(t_1-t_0)\otimes\dots\otimes\Cal E_P(t_n-t_{n-1})$
defined by 
$$
u_i = T_i(t_1-t_0)\otimes\dots\otimes T_i(t_n-t_{n-1}), 
$$
$i=1,2$.  Notice that because of (2.8) we have 
$\|u_i\|\leq 1$ for $i=1,2$, and  
$$
f_\Cal P(T_i,T_j;t) = 
\<u_i,u_j\>.  
$$
By an obvious induction using nothing more than 
the associative law, Theorem 1.12 implies that 
there is a unique multiplication operator 
$$
M:\Cal E_P(t_1-t_0)\otimes\dots\otimes\Cal E_P(t_n-t_{n-1}) \to
\Cal E_P(t)
$$
satisfying $M(v_1\otimes\dots\otimes v_n) = v_1v_2\dots v_n$, 
and moreover that $\|M\|\leq 1$.  Noting that 
$Mu_i = T_i(t)$ and using $\|M\|\leq 1$ we have 
$$
\align
&|f_\Cal P(T_i,T_j;t)-\<T_i(t),T_j(t)\>| = 
|\<u_i,u_j\> - \<Mu_i,Mu_j\>| \\
&= 
|\<(\bold 1-M^*M)u_i,u_j\>| 
\leq \|(\bold 1-M^*M)^{1/2}u_i\| \cdot
\|(\bold 1-M^*M)^{1/2}u_j\|.  
\endalign
$$
Since 
$$
\align
\|(\bold 1-M^*M)^{1/2}u_j\|^2 &= \<(\bold 1-M^*M)uj, u_j\> = 
\|u_j\|^2 - \|Mu_j\|^2  \\
&\leq 1-\|Mu_j\|^2= 1-\<T_j(t),T_j(t)\>,   
\endalign
$$
the estimate of Lemma 2.14 follows \qed
\enddemo

\demo{proof of Lemma 2.15}
We show first that for every unit $T\in \Cal U_P$, 
$$
\lim_{t\to0+}\<T(t),T(t)\> = 1.  \tag{2.17}
$$
Indeed, since units must satisfy a growth condition of the 
form $\<T(t),T(t)\>\leq e^{kt}$ it suffices to show that
$$
1\leq \liminf_{t\to 0+}\<T(t),T(t)\>.  \tag{2.18}
$$
Now for every $t>0$ the map
$$
x\in \Cal B(H)\mapsto \<T(t),T(t)\>P_t(x) - T(t)xT(t)^*
$$
is completely positive; taking $x=\bold 1$ we find that
for every unit vector $\xi\in H$ 
$$
\|T(t)^*\xi\|^2=\<T(t)T(t)^*\xi,\xi\>
\leq \<T(t),T(t)\>\<P_t(\bold 1)\xi,\xi\>.  
$$
As $t\to 0+$, $\<P_t(\bold 1)\xi,\xi\>$ tends to 
$\<\bold 1\xi,\xi\> = 1$, and since $T(t)^*\xi$ tends 
to $\xi$ in the norm of $H$ we have $\|T(t)^*\xi\|\to 1$.  
(2.18) follows.  

Now let $T_1,T_2\in \Cal U_P$.  
Because each unit satisfies a growth condition of the 
form (2.3) and since we can replace each $T_j(t)$ by 
$e^{-k_jt}T_j(t)$ without affecting the conclusion of Lemma
2.15, it suffices to prove Lemma 2.15 for units 
$T_1,T_2$ satisfying $\<T_j(t),T_j(t)\>\leq 1$ for all $t>0$.  
Fix such a pair $T_1, T_2$, fix $t>0$, and set 
$$
u=T_1(t), \quad v = 
\<T_1(t),T_1(t)\>T_2(t)-\<T_2(t),T_1(t)\>T_1(t).  \tag{2.19}
$$
$u$ and $v$ are orthogonal elements of $\Cal E_P(t)$.  
We claim that for any two orthogonal elements $u,v\in \Cal E_P(t)$ 
we have 
$$
\<u,u\>vv^*\leq \<v,v\>(\<u,u\>P_t(\bold 1)-uu^*).  \tag{2.20}
$$
Indeed, (2.20) is trivial if either $u$ or $v$ is $0$, so we 
assume that both are nonzero.  In this case, put 
$$
u_0 = \<u,u\>^{-1/2}u, \quad v_0 = \<v,v\>^{-1/2}v.  
$$
Then $\{u_0, v_0\}$ is part of an orthonormal basis for 
$\Cal E_P(t)$, hence the map
$$
x\mapsto P_t(x) - u_0xu_0^* - v_0xv_0^*
$$
is completely positive.  Taking $x=\bold 1$ we find that 
$$
v_0v_0^* \leq P_t(\bold 1) - v_0v_0^*,
$$
and (2.20) follows after multiplying through by 
$\<u,u\>\<v,v\>$.  

For $u$ and $v$ as in (2.19), the inequality (2.20) implies 
that for every unit vector $\xi\in H$, 
$$
\align
\<T_1(t),T_1(t)\>&\|\<T_1(t),T_1(t)\>T_2(t)^*\xi - 
\<T_1(t),T_2(t)\>T_1(t)^*\xi\|^2 \leq\\ 
&\<v,v\>(\<T_1(t),T_1(t)\>\<P_t(\bold 1)\xi,\xi\>-\|T_1(t)^*\xi\|^2).  
\endalign
$$
Notice that $\<v,v\>\leq 4$.  Indeed, since 
$\|T_j(t)\|_{\Cal E_P(t)}\leq \<T_j(t),T_j(t)\>^{1/2}\leq 1$ 
we have 
$$
\<v,v\>= 
\|\<T_1(t),T_1(t)\>T_2(t)-\<T_2(t),T_1(t)\>T_1(t)\|^2_{\Cal E_P(t)} 
\leq 4.  
$$
Thus the preceding inequality implies that 
$$
\|\<T_1(t),T_1(t)\>T_2(t)^*\xi - 
\<T_1(t),T_2(t)\>T_1(t)^*\xi\|^2   \tag{2.21}
$$
is dominated by a term of the form 
$$
\frac{4}{\<T_1(t),T_1(t)\>}
(\<T_1(t),T_1(t)\>\<P_t(\bold 1)\xi,\xi\>-\|T_1(t)^*\xi\|^2). \tag{2.22}
$$
As $t\to 0+$, the expression in (2.22) tends to zero because of 
(2.17) and the fact that both $\<P_t(\bold 1)\xi,\xi\>$ and 
$\|T_1(t)^*\xi\|^2$ tend to $\|\xi\|^2 = 1$.  Thus the term in 
(2.21) tends to zero as $t\to 0+$.  Taking note of (2.17) once again, 
we conclude that 
$$
\lim_{t\to 0+}\|T_2(t)^*\xi - \<T_1(t),T_2(t)\>T_1(t)^*\xi\| = 0.  
$$
Writing 
$$
\align
|1-\<T_1(t),T_2(t)\>| = &\|\xi -\<T_1(t),T_2(t)\>\xi\| \leq \\
&\|\xi - T_2(t)^*\xi\| + |\<T_1(t),T_2(t)\>|\cdot\|\xi-T_1(t)^*\xi\| +\\
&\|T_2(t)^*\xi - \<T_1(t),T_2(t)\>T_1(t)^*\xi\|,    
\endalign
$$
and noting that each of the three terms on 
the right tends to zero as $t\to 0+$, we obtain
$$
\lim_{t\to 0+}|1 - \<T_1(t),T_2(t)\>| = 0
$$
as required for Lemma 2.15\qed 
\enddemo
That also completes the proof of Theorem 2.6 \qed
\enddemo

\subheading{3.  Lifting units}
Let $\alpha = \{\alpha_t: t\geq 0\}$ be an \esg\ acting on 
$M = \Cal B(H)$, $H$ separable.  
$\alpha$ can be compressed to {\it certain} hereditary 
subalgebras $M_0 = p_0Mp_0$ of $M$ so as to give a CP semigroup
$P$ acting on $M_0 \cong \Cal B(p_0H)$.  In this section we
show that the units of $\alpha$ map naturally to those of 
$P$, and in the case where $\alpha$ is minimal over $P$ we
show that this map is a bijection (Theorem 3.6).  

A projection $p_0\in M$ is said to be {\bf increasing} if 
$\alpha_t(p_0) \geq p_0$ for every $t\geq 0$.  In this case
we obtain a CP semigroup $P = \{P_t: t\geq 0\}$ acting on 
$M_0 = p_0Mp_0$ by way of 
$$
P_t(x) = p_0\alpha_t(x)p_0, \qquad t\geq 0, x\in M_0.  
$$
$P$ is called a {\bf compression} of $\alpha$ and $\alpha$ 
is called a {\bf dilation} of $P$.  It is possible for $P$ 
itself to be an \esg , that is to say 
$P_t(xy) = P_t(x)P_t(y)$ for every $x,y\in M_0$, $t\geq 0$.  
In this case we call $P$ a multiplicative compression of 
$\alpha$.  Finally, $\alpha$ is said to be {\bf minimal} over
$P$ if there are no intermediate multiplicative compressions;
more explicitly, there should exist no
increasing projection $q\in M$ for which a) $q\geq p_0$ and 
b) the compression of $\alpha$ to $qMq$ is multiplicative, other
than $q=\bold 1$.  

The issue of minimality over $P$ merits some discussion (for 
full details see \cite{3}).  The condition
$$
\alpha_t(p_0) \uparrow \bold 1_H, \quad \text{ as }t \to \infty 
$$
is necessary, but {\it not} sufficient for minimality.  There are
a number of equivalent additional conditions that guarantee 
minimality, and the one we require is formulated as follows.  
For every $t>0$, let $q_t$ be the projection onto the subspace
$[\alpha_t(M)p_0H]$.  $q_t$ obviously belongs to the commutant 
of $\alpha_t(M)$.  For every 
fixed $t>0$ and every partition $\Cal P= \{0=t_0<t_1<\dots<t_n=t\}$
of the interval $[0,t]$, we set 
$$
q_{\Cal P, t} = 
q_{t_1}\alpha_{t_1}(q_{t_2-t_1})\alpha_{t_2}(q_{t_3-t_2})
\dots\alpha_{t_{n-1}}(q_{t_n-t_{n-1}}).  \tag{3.1}
$$
It is shown in \cite{3, Proposition 3.4} that $q_{\Cal P, t}$ is a
projection
in the commutant of $\alpha_t(M)$ and that 
$$
\Cal P_1\subseteq \Cal P_2 \implies 
q_{\Cal P_1, t}\leq q_{\Cal P_2, t}.
$$
Thus the strong limit 
$$
\bar q_t = \lim_{t=\to\infty}q_{\Cal P, t}
$$
exists for every $t>0$ and the resulting family
of projections $\{\bar q_t\in \alpha_t(M)^\prime: t>0\}$ satisfies
the cocycle equation
$$
\bar q_{s+t} = \bar q_s\alpha_t(\bar q_t), \qquad s, t >0
$$
as well as a natural continuity condition.  Moreover, it was
shown in \cite{3} that $\alpha$ is minimal over 
$P$ iff the following two conditions are satisfied
$$
\align
\alpha_t(p_0) \uparrow &\bold 1,\qquad \text{as } t\to \infty, \tag{3.2.1}\\
\bar q_t = &\bold 1, \qquad \text{for every }t > 0.  \tag{3.2.2}
\endalign
$$

The purpose of this section is to show how the units of $\alpha$ 
are related to the units of $P$ in the case where $\alpha$ is 
minimal over $P$.  More precisely, let
$\Cal E_\alpha = \{\Cal E_\alpha(t) : t>0\}$ be the product system 
of $\alpha$.  Thus $\Cal E_\alpha(t)$ is the intertwining space
$$
\Cal E_\alpha(t) = \{T\in \Cal B(H): \alpha_t(x)T = Tx, x\in \Cal B(H)\}
$$
which becomes a separable Hilbert space with respect to the
inner product defined by 
$$
\<S, T\> \bold 1 = T^*S, \qquad S, T\in \Cal E_\alpha(t).  
$$

\proclaim{Proposition 3.3}
For every $t>0$ and every operator $T\in \Cal E_\alpha(t)$, the subspace 
$p_0H$ is invariant under the adjoint $T^*$.  The operator 
$S\in \Cal B(p_0H)$ defined by 
$$
S^* = T^*\restriction_{p_0H}
$$
belongs to the space $\Cal E_P(t)$ and satisfies
$$
\<S, S\>_{\Cal E_P(t)} \leq \<T, T\>_{\Cal E_\alpha(t)}.  \tag{3.4}
$$
\endproclaim
\demo{proof}
The proof is a straightforward consequence of the fact that
$p_0$ is an increasing projection.  Indeed, if we choose an 
orthonormal basis $\{v_1, v_2, \dots\}$ for $\Cal E_\alpha(t)$ then
we have 
$$
\sum_k v_k(\bold 1- p_0) v_k^*  = \alpha_t(\bold 1-p_0) 
= \bold 1 - \alpha_t(p_0) \leq \bold 1-p_0.  
$$
It follows that $v_k(\bold 1-p_0)p_k^*\leq \bold 1-p_0$ for every $k$, 
hence $v_k$ leaves the orthogonal complement of $p_0H$ invariant 
for every $k$, and hence $v_k^*p_0H\subseteq p_0H$.  
Since the linear span of the $\{v_k\}$ 
is dense in $\Cal E_\alpha(t)$ in the operator norm, the 
assertion $\Cal E_\alpha(t)^* p_0H \subseteq p_0H$ follows.  

Let $S$ be the indicated operator in $\Cal B(p_0H)$.  Since 
$\alpha$ is an \esg\ the Hilbert space norm of an element of 
$\Cal E_\alpha(t)$ coincides with its operator norm.  Thus,
in order to 
show that $S\in \Cal E_P(t)$ and satisfies the inequality (3.4), 
it suffices to show that the operator mapping $L$ of $\Cal B(p_0H)$ defined
by 
$$
L(x) = \|T\|^2 P(x) - SxS^*
$$
is completely positive.  Now by definition of $S$ see that 
for every $x\in p_0Mp_0$ we have 
$$
SxS^* = p_0TxT^*p_0 = p_0\alpha_t(x)TT^*p_0.  
$$
Since $TT^*$ is a positive operator of norm $\|T\|^2$ in the
commutant of $\alpha_t(M)$ it follows that $C = (\|T\|^2\bold 1 - TT^*)^{1/2}$
is a positive operator in the commutant of $\alpha_t(M)$, hence 
$$
L(x) = \|T\|^2 p_0\alpha_t(x)p_0 - p_0\alpha_t(x)TT^*p_0
= p_0C\alpha_t(x) Cp_0
$$
is obviously a completely positive mapping of $p_0Mp_0$ into itself \qed
\enddemo

Proposition (3.3) implies that there is a natural mapping of the units 
of $\alpha$ to the units of $P$, defined as follows.  In this concrete 
setting we may consider a unit of $\alpha$ to be a strongly continuous
semigroup $T = \{T(t): t\geq 0\}$ of operators in $M$ satisfying 
$$
\alpha_t(x)T(t) = T(t)x, \qquad x\in M.  
$$
Choose such a $T$, and for every $t>0$ let 
$S(t)\in \Cal B(p_0H)$ be the operator defined by
$$
S(t)^* = T(t)^*\restriction p_0H.  \tag{3.5}
$$
It is obvious that $S = \{S(t) : t>0\}$ is a strongly continuous
semigroup of bounded operators on $p_0H$ for which $S(t)\to \bold 1$ 
strongly as $t\to 0+$, and we have $S(t) = p_0T(t)p_0 = p_0T(t)$ for
every $t$.  Proposition (3.3) implies that $S(t)$ 
belongs to $\Cal E_P(t)$ for every $t>0$ and moreover
$$
\<S(t), S(t)\>_{\Cal E_P(t)} \leq \|T(t)\|^2.  
$$
Because $T$ is a unit of $\alpha$ we must have 
$$
\|T(t)\|^2 = e^{tC(T,T)}, \qquad t>0
$$
where $C: \Cal U_\alpha \times \Cal U_\alpha \to \Bbb C$ is the 
covariance function of $\alpha$, and thus $S$ is a unit of $P$.  

\proclaim{Theorem 3.6}
Suppose that $\alpha$ is minimal over $P$.  Then the function 
$\theta: \Cal U_\alpha \to \Cal U_P$ defined by $\theta(T) = S$
is a bijection.  
\endproclaim

\demo{proof}
In order to show that $\theta$ is one-to-one, 
fix $T_1, T_2\in \Cal U_\alpha$ such that $\theta(T_1) = \theta(T_2)$.  
Thus $T_1(t)^*\restriction_{p_0H} = T_2(t)^*\restriction_{p_0H}$,
for every $t>0$.  Noting that $\alpha_t(x)T_k(t) = T_k(t)x$ it 
follows that for every $x\in M$ and $\xi\in p_0H$ we have 
$$
T_1^*(t)\alpha_t(x)\xi = xT_1(t)^*\xi = xT_2(t)^*\xi = 
T_2^*(t)\alpha_t(x)\xi.  
$$
Letting $q_t$ be the projection on the subspace $[\alpha_t(M)p_0H]$ and
taking adjoints, the preceding formula implies that
$$
q_tT_1(t) = q_tT_2(t), \qquad t>0.  
$$

Note too that the preceding formula implies that for every $0<s<t$
we have 
$$
q_s\alpha_s(q_{t-s})T_1(t) = q_s\alpha_s(q_{t-s})T_2(t).  \tag{3.7}
$$
Indeed, the left side of (3.7) can be written
$$
\align
q_s\alpha_s(q_{t-s})T_1(s)T_1(t-s) &= q_sT_1(s)q_{t-s}T_1(t-s) 
= q_sT_2(s)q_{t-s}T_2(t-s) \\
&= q_s\alpha_s(q_{t-s})T_2(s)T_2(t-s)
\endalign
$$
and (3.7) follows.  By an obvious induction argument, it follows
similarly that if $\Cal P = \{0 = t_0 < t_1 < \dots < t_n = t\}$ 
is any finite partition of the interval $[0, t]$ and if 
$q_{\Cal P, t}$ is defined as in the discussion above, then we have 
$$
q_{\Cal P, t}T_1(t) = q_{\Cal P, t}T_2(t).  
$$
Because of the minimality condition (3.2.2) we may take the 
limit on $\Cal P$ to obtain $T_1(t) = T_2(t)$.  

In order to show that $\theta$ is surjective, we require the following

\proclaim{Lemma 3.8}
Let $S = \{S(t): t\geq 0\}$ be a unit of $P$ and for every $t>0$ 
let $q_t$ be the projection onto $[\alpha_t(M)p_0H]$.  

Then for every $t>0$ there is a unique operator $v_t\in \Cal E_\alpha(t)$
satisfying the two conditions $q_tv_t = v_t$, and 
$v_t^*\restriction p_0H = S_t^*$.  Moreover, there is a real 
constant $k$ such that $\|v_t\| \leq e^{kt}$ for every $t>0$.  
\endproclaim

\demo{proof}
Let $S=\{S(t): t\geq 0\}$ be a semigroup of bounded operators on
$\Cal B(p_0H)$ and let $k$ be a real number
with the property that for every $t>0$, 
$$
z\in \Cal B(p_0H) \mapsto e^{kt}P_t(z)-S(t)zS(t)^* \tag{3.9}
$$
is a completely positive map.   
Let $x_1, x_2, \dots, x_n$ be a set of operators in the
larger von Neumann algebra $M=\Cal B(H)$ and choose vectors
$\xi_1, \xi_2,\dots,\xi_n\in p_0H$.  We claim
$$
\|\sum_{k=1}^nx_kS(t)^*\xi_k\|^2 \leq 
e^{kt}\|\sum_{k=1}^n \alpha_t(x_k)\xi_k\|^2. \tag{3.10}
$$
Indeed, the left side of (3.10) is 
$$
\sum_{k,j=1}^n\<x_kS(t)^*\xi_k, x_jS(t)^*\xi_j\> =
\sum_{k,j=1}^n\<S(t)p_0x_j^*x_kp_0S(t)^*\xi_k,\xi_j\>.  \tag{3.11}
$$
Since the $n\times n$ matrix $(a_{jk})$ defined by 
$a_{jk} = p_0x_j^*x_kp_0$ is a positive operator matrix with
entries from $p_0Mp_0$, (3.9) implies that 
the right side of (3.11) is dominated by 
$$
e^{kt}\sum_{k,j=1}^n\<\alpha_t(p_0x_j^*x_kp_0)\xi_k, \xi_j\> =
e^{kt}\|\sum_{k=1}^n \alpha_t(x_kp_0)\xi_k\|^2.  
$$
Since $p_0$ is an increasing projection and $\xi_k\in p_0H$, we
can write $\alpha_t(x_kp_0)\xi_k = \alpha_t(x_k)\alpha_t(p_0)\xi_k 
= \alpha_t(x_k)\xi_k$ for each $k=1,2,\dots n$, and hence the 
right side of the previous formula becomes
$$
e^{kt}\|\sum_{k=1}^n \alpha_t(x_k)\xi_k\|^2.  
$$
The inequality (3.10) follows.  

From (3.10) it follows that there is a unique operator 
$v_t\in \Cal B(H)$, having norm at most $e^{kt/2}$, and which 
satisfies 
$$
\align
v_t^*\alpha_t(x)\xi &= xS(t)^*\xi, \qquad x\in \Cal B(H), 
\xi\in p_0H, \text{ and}  \tag{3.12.1}\\
v_t^* &= v_t^*q_t.  \tag{3.12.2}
\endalign
$$

We claim that $v_t\in \Cal E_\alpha(t)$ or equivalently, that
$$
v_t^*\alpha_t(x) = xv_t^*, \qquad x\in \Cal B(H).  \tag{3.13}
$$
Indeed, because of (3.12.2) we have 
$v_t^*\alpha_t(x) = v_t^*q_t\alpha_t(x) = v_t^*\alpha_t(x)q_t$, 
and similarly $xv_t^* = xv_t^*q_t$.  Thus it suffices to show that
the operators on both sides of (3.13) agree on vectors in 
$q_tH = [\alpha_t(M)p_0H]$.  If such a vector has the form
$\eta = \alpha_t(y)\xi$ with $y\in \Cal B(H)$ and $\xi\in p_0H$ then
we have
$$
v_t^*\alpha_t(x)\eta = v_t^*\alpha_t(x)\alpha_t(y)\xi = 
v_t^*\alpha_t(xy)\xi = xyS(t)^*\xi = xv_t^*\alpha_t(y)\xi, 
$$
and (3.13) follows because such vectors $\eta$ span the range of $q_t$.  

This proves the existence assertion of Lemma 3.8.  For uniqueness, let
$w_t\in \Cal E_\alpha(t)$ satisfy $q_tw_t = w_t$ and 
$w_t^*\restriction_{p_0H} = S(t)^*$.  Then for any vector $\eta$ of the
form $\eta = \alpha_t(x)\xi$, $x\in \Cal B(H)$, $\xi\in p_0H$ we have 
$$
w_t^*\eta = w_t^*\alpha_t(x)\xi = xw_t^*\xi = xS(t)^*\xi =
v_t^*\eta,  
$$
so that $w_t^*$ and $v_t^*$ agree on $[\alpha_t(M)p_0H] = q_tH$, and
hence $w_t^* = w_t^*q_t$ and $v_t^* = v_t^*q_t$ agree \qed 
\enddemo

To complete the proof of Theorem 3.6, choose a unit $S = \{S_t: t\geq 0\}$ 
for $P$ and let $\{v_t: t>0\}$ be the family of operators defined by 
Lemma 3.8. This family of operators is certainly a section of the product
system of $\alpha$, but it is not a unit because 
it does not satisfy the semigroup property 
$v_{s+t} = v_sv_t$.  In order to obtain a unit from this family 
$\{v_t: t>0\}$ we carry out the following construction.  

Fix $t>0$.  For every finite partition 
$\Cal P = \{0 = t_0 < t_1 < \dots < t_n = t\}$ of the interval $[0, t]$, 
consider the operator
$$
v_{\Cal P, t} = v_{t_1-t_0}v_{t_2-t_1}\dots v_{t_n-t_{n-1}}.  
$$
It is clear that $v_{\Cal P,t}$ belongs to $\Cal E_\alpha(t)$, and 
because of the growth condition $\|v_s\| \leq e^{ks}$ for all positive
$s$ we have 
$$
\|v_{\Cal P,t}\| \leq e^{kt}.  
$$
Thus $\Cal P \mapsto v_{\Cal P,t}$ defines a bounded net of operators
belonging to the $\sigma$-weakly closed operator space $\Cal E_\alpha(t)$.  
We will show next that this net converges weakly.  The resulting limit 
$$
T_t = \lim_{\Cal P}v_{\Cal P,t}
$$
will satisfy the semigroup property $T_{s+t} = T_sT_t$, but since
the net of finite partitions is uncountable, 
continuity (or even measurability) in $t$ is not 
immediate.  We then give a separate argument which guarantees that 
$\{T_t: t>0\}$ is strongly continuous, and that the unit of $\alpha$ that
it defines maps to $S$ as required.  

\proclaim{Lemma 3.14}
For every $t>0$ and every finite partition $\Cal P$ of $[0,t]$, let
$q_{\Cal P, t}$ be the projection defined in (3.1).  Then for 
every pair of partitions satisfying $\Cal P_1 \subseteq \Cal P_2$ we
have 
$$
q_{\Cal P_1}v_{\Cal P_2, t} = v_{\Cal P_1,t}.  
$$
\endproclaim

\remark{Remark 3.15}
We have already seen that the net of 
projections $\Cal P \mapsto q_{\Cal P,t}$ 
is increasing in $\Cal P$ and by minimality of $\alpha$ over $P$ 
this net of projections has limit $\bold 1$ for every fixed $t>0$.  
Thus the coherence condition asserted in  
Lemma 3.14, together with the fact that 
$\|v_{\Cal P,t}\| \leq e^{kt}$, implies that the net of 
adjoint operators 
$$
\Cal P \mapsto (v_{\Cal P,t})^*
$$
must converge in the {\it strong} operator topology.  In particular,
the weak limit
$$
T_t = \lim_\Cal P v_{\Cal P, t}
$$
exists for every $t$ and defines an element of $\Cal E_\alpha(t)$.  
\endremark

\demo{proof of Lemma 3.14}
We claim first that for every $s$, $t>0$ we have
$$
q_{s+t}v_sv_t = v_{s+t}.  \tag{3.16}
$$
Indeed, because of the uniqueness assertion of Lemma (3.8), it 
suffices to show that the operator $w= q_{s+t}v_sv_t$ belongs 
to $\Cal E_\alpha(s+t)$ and satisfies $w^*\restriction_{p_0H} = S(s+t)^*$.  
The first assertion is obvious because $v_sv_t\in \Cal E_\alpha(s+t)$ 
and $q_{s+t}$ commutes with $\alpha_{s+t}(M)$.  To see that $w^*$ 
restricts to $S(s+t)^*$ on $p_0H$, choose $\xi \in p_0H$ and note that
$$
w^*\xi = v_t^*v_s^*q_{s+t}\xi = v_t^*v_s^*\xi = v_t^*S(s)^*\xi =
S(t)^*S(s)^*\xi = S(s+t)^*\xi.  
$$
Thus (3.16) is established.  

In order to prove Lemma (3.14), it is enough to consider the 
case where $\Cal P_2$ is obtained from  
$\Cal P_1 = \{0=t_0 < t_1 < \dots < t_n = t\}$ by adjoining to it 
a single point $\tau$, say 
$$
t_k < \tau < t_{k+1}
$$ 
for some $k = 0, 1, \dots, n-1$.  Now by (3.16) we see that
$$
q_{t_{k+1}-t_k}v_{\tau - t_k}v_{t_{k+1}-\tau} = v_{t_{k+1}-t_k},
$$
and if we make this substitution for $v_{t_{k+1}-t_k}$ in the formula
$$
v_{\Cal P_1,t} = v_{t_1-t_0}\dots 
v_{t_{k+1}-t_k}\dots v_{t_n-t_n-1} 
$$
we obtain 
$$
\align
v_{\Cal P_1,t}
&= 
v_{t_1-t_0}\dots v_{t_k-t_{k-1}}(q_{t_{k+1}-t_k}v_{\tau - t_k}v_{t_{k+1-\tau}})
v_{t_{k+2}-t_{k+1}}\dots v_{t_n-t_{n-1}} \\
&= 
(q_{t_1-t_0}v_{t_1-t_0})\dots 
(q_{t_{k+1}-t_k}v_{\tau-t_k}v_{t_{k+1}-\tau})
\dots (q_{t_n-t_{n-1}}v_{t_n-t_{n-1}}).  
\endalign
$$
If we now move each of the ``$q$" terms to the left, using the relation
$v_sx = \alpha_s(x)v_s$, $x\in \Cal B(H)$, that last expression on the 
right becomes
$$
q_{t_1-t_0}\alpha_{t_1}(q_{t_2-t_1})\dots
\alpha_{t_{n-1}}(q_{t_n-t_{n-1}})v_{t_1-t_0}\dots v_{\tau-t_k}v_{t_{k+1-\tau}}
\dots v_{t_n-t_{n-1}},
$$
which is $q_{\Cal P_1,t}v_{\Cal P_2,t}$, as required in Lemma 3.14 \qed
\enddemo

It follows from Remark 3.15 that we have {\it strong} convergence of the 
net of adjoints 
$$
T_t^* = \lim_\Cal P(v_{\Cal P,t})^* 
$$
for every positive $t$.  Since multiplication is strongly continuous
on bounded sets we obtain $T_t^*T_s^*$ as a strong double limit
$$
T_t^*T_s^* = \lim_{\Cal P_1, \Cal P_2} (v_{\Cal P_1,t})^*(v_{\Cal P_2,s})^* =
\lim_{\Cal P_1, \Cal P_2} (v_{\Cal P_1,s}v_{\Cal P_2,t})^*.  
$$
Taking adjoints, we have the following weak convergence
$$
T_sT_t = \lim_{\Cal P_1, \Cal P_2}v_{\Cal P_1,s}v_{\Cal P_2,t}
=  \lim_{\Cal P_1, \Cal P_2}v_{\Cal P_1\cup (s+\Cal P_2),s+t},   
$$
where $\Cal P_1\cup (s+\Cal P_2)$ denotes the partition of $[0,s+t]$ 
obtained by first listing the elements of $\Cal P_1$ and then 
listing the elements of $s+\Cal P_2$.  
Since the right side is a limit over a cofinal subnet of partitions of 
the interval $[0, s+t]$, we conclude that $T_sT_t = T_{s+t}$ for every
positive $s,t$.  

We claim next that $T_t^*\restriction_{p_0H} = S(t)^*$.  To see this, 
notice that since $v_s^*$ restricts to $S(s)^*$ for every positive $s$ 
and $\{S(s) : s\geq 0\}$ is a semigroup, it follows that 
$(v_{\Cal P, t})^*$ restricts to $S(t)^*$ for every $t>0$.  The claim 
follows because the net $(v_{\Cal P, t})^*$ converges weakly to 
$T_t^*$.  

Finally, we show that the semigroup $\{T_t^*: t>0\}$ is strongly 
continuous; that is, we will show that
$$
\lim\|T_t^*\xi - \xi\| = 0, \tag{3.17}
$$
for every $\xi \in H$.  Indeed, (3.17) is certainly true in case 
$\xi\in p_0H$, because $T_t^*$ restricts to $S(t)^*$ and $S$ is 
a continuous semigroup of operators on $p_0H$.  Let $K$ denote the 
set of all vectors $\xi\in H$ for which (3.17) holds.  $K$ is clearly
a closed subspace of $H$ which contains $p_0H$.  We assert now that
for every $s>0$, 
$$
\alpha_s(M)K\subseteq K. \tag{3.18}
$$
Indeed, if $s>0$ and $x\in M=\Cal B(H)$, then for sufficiently small
positive $t$ we have $t<s$ and hence
$$
T_t^*\alpha_s(x)\xi = \alpha_{s-t}(x)v_t^*\xi.  
$$
So if $\xi\in K$ then 
$$
\align
\|T_t^*\alpha_s(x)\xi - \alpha_s(x)\xi\| &=
\|\alpha_{s-t}(x)v_t^*\xi - \alpha_s(x)\xi\|  \\
&\leq \|\alpha_{s-t}(x) v_t^*\xi - \alpha_{s-t}(x)\xi\| +
\|\alpha_{s-t}(x)\xi - \alpha_s(x)\xi\|  \\
&\leq \|x\|\cdot\|v_t^*\xi - \xi\| + \|\alpha_{s-t}(x)\xi - \alpha_s(x)\xi\|. 
\endalign 
$$
Both terms on the right tend to $0$ with $t$ because $\xi\in K$ and 
$\alpha$ is a (continuous) \esg. Thus $K$ contains every vector of 
the form $\alpha_s(x)p_0\xi$, where $x\in \Cal B(H)$ and $\xi\in H$ 
are arbitrary, and $s$ is an arbitrary positive number.  Allowing 
$s$ to tend to zero we find that $\alpha_s(x)$ tends strongly to $x$, 
and hence
$$
K \supseteq [\Cal B(H)p_0H] = H.  
$$
Thus $\{T_t: t>0\}$ is strongly continuous.  

It follows that $u = \{T_t: t>0\}$ is a unit of $\alpha$ for which
$\theta(u) = S$, and the proof of Theorem 3.6 is complete.  
\enddemo

\remark{Remark 3.18}
Notice that the semigroup $T = \{T_t: t > 0\} \in \Cal U_\alpha$
defined by 
$$
T_t = \lim_{\Cal P}v_{\Cal P,t}, \qquad t>0
$$
projects as follows relative to any finite partition 
$\Cal P = \{0=t_0 < t_1 < \dots < t_n=t\}$ of $[0,t]$: 
$$
q_{\Cal P,t}T_t = v_{\Cal P,t} = 
v_{t_1-t_0}v_{t_2-t_1}\dots v_{t_n-t_{n-1}}.  \tag{3.19}
$$
\endremark

\subheading{4.  The covariance function of a CP semigroup}
Let $P = \{P_t: t\geq 0\}$ be a {\it unital} CP semigroup acting on 
$\Cal B(H_0)$.  By a theorem of B. V. R. Bhat, there is a Hilbert 
space $H$ containing $H_0$ and an \esg\ 
$\alpha = \{\alpha_t: t\geq 0\}$ acting on 
$\Cal B(H)$ such that the projection $p_0$ onto
$H_0$ is increasing for $\alpha$ and $P$ is obtained by compressing
$\alpha$ to $p_0\Cal B(H)p_0 \cong \Cal B(p_0H)$ as we have described
above \cite{6.7}.  Moreover, one may also arrange 
(by passing to a suitable intermediate
\esg\ if necessary) that $\alpha$ is minimal over $P$ \cite{3}.  
Finally, any two minimal dilations of $P$ are conjugate.  

The purpose of this section is to calculate the covariance function
$$
c_P: \Cal U_P\times \Cal U_P \to \Bbb C
$$
of $P$ in terms of the covariance function 
$$
c_\alpha: \Cal U_\alpha \times \Cal U_\alpha \to \Bbb C
$$ of $\alpha$ when $\alpha$ is the minimal dilation of $P$.  Indeed,
letting $\theta: \Cal U_\alpha \to \Cal U_P$ be the bijection defined
by Theorem 3.6, we will show that 
$$
c_P(\theta(u_1), \theta(u_2)) = c_\alpha(u_1, u_2).  \tag{4.1}
$$
Once one has (4.1), it is apparent that the bijection 
$\theta$ gives rise to a natural unitary operator from the Hilbert 
space associated with $(\Cal U_\alpha, c_\alpha)$ onto 
that associated with $(\Cal U_P, c_P)$, and in particular, these
two Hilbert spaces have the same dimension.  Hence, 
{\it the numerical index $d_*(P)$ of $P$ must agree with the numerical
index $d_*(\alpha)$ of its minimal dilation $\alpha$}.  

For every $t>0$ and every partition 
$\Cal P = \{0=t_0 < t_1 < \dots < t_n\}$ let $q_{\Cal P,t}$ 
be the projection defined by (3.1).  Since $q_{\Cal P,t}$ 
belongs to the commutant of $\alpha_t(\Cal B(H))$ it follows 
that 
$$
q_{\Cal P,t}\Cal E_\alpha(t) \subseteq \Cal E_\alpha(t).  
$$
Thus we may consider the left multiplication operator 
$$
Q_{\Cal P,t}: x\in \Cal E_\alpha(t) 
\mapsto q_{\Cal P,t}x \in \Cal E_\alpha(t)
$$
as a bounded operator on the Hilbert space $\Cal E_\alpha(t)$.  
$Q_{\Cal P,t}$ is a self-adjoint projection in 
$\Cal B(\Cal E_\alpha(t))$.  

\proclaim{Proposition 4.2}
The projections $Q_{\Cal P,t}\in \Cal B(\Cal E_\alpha(t))$ are 
increasing in the variable $\Cal P$ and 
$$
\lim_{\Cal P}\|Q_{\Cal P,t}x - x\|_{\Cal E_\alpha(t)} = 0, 
\qquad x\in \Cal E_\alpha(t).  
$$
\endproclaim
\demo{proof}
These assertions are a simple consequence of the definition of the 
inner product $\<\cdot, \cdot \>$ in $\Cal E_\alpha(t)$:
$$
\<S,T\> \bold 1 = T^*S, \qquad S,T\in \Cal E_\alpha(t).  
$$
Indeed, if $\Cal P_1$ and $\Cal P_2$ are two finite partitions 
of $[0,t]$ satisfying $\Cal P_1 \subseteq \Cal P_2$, then for every
operator $T \in \Cal E_\alpha(t)$ we have 
$$
\<Q_{\Cal P_1,t}T,T\> \bold 1_H = 
T^*q_{\Cal P_1,t}T \leq T^*q_{\Cal P_2,t}T = 
\<Q_{\Cal P_2,t}T, T\>\bold 1_H,  
$$
hence $Q_{\Cal P_1,t} \leq Q_{\Cal P_2,t}$.  Similarly, the
fact that the net 
$\Cal P \mapsto Q_{\Cal P,t} \in \Cal B(\Cal E_\alpha(t))$ converges
to the identity of $\Cal B(\Cal E_\alpha(t))$ follows immediately
from (3.2.2) \qed
\enddemo

In section 2, the covariance function of a CP semigroup $P$ is 
defined in terms of limits of certain finite products of complex
numbers of the form 
$$
\<S_1(t),S_2(t)\>_{\Cal E_P(t)} = \<S_1(t),S_2(t)\>.$$
We now 
show how these products are expressed in terms of $\alpha$.  

\proclaim{Theorem 4.3}
Let $S_1$ and $S_2$ be two units of a unital CP semigroup
$P$.  Let $\alpha$ be its minimal dilation to an \esg\ and let 
$T_1$, $T_2$ be the unique units of $\alpha$ satisfying 
$\theta(T_k) = S_k$, $k=1,2$. 

Then for every $t>0$ and every finite partition 
$\Cal P = \{0=t_0 < t_1 < \dots < t_n=t\}$ 
of the interval $[0,t]$, we have 
$$
\prod_{k=1}^n\<S_1(t_k-t_{k-1}),S_2(t_k-t_{k-1})\>
=\<Q_{\Cal P,t}T_1(t), T_2(t)\>,  
$$
the inner product on the right being relative to 
the Hilbert space $\Cal E_\alpha(t)$.  
\endproclaim

\demo{proof}For each $t>$, let $q_t$ be the projection
onto the subspace $[\alpha_t(\Cal B(H)p_0H]$.  
Lemma 3.8 guarantees that
there is a unique pair of operators 
$v_1(t), v_2(t) \in \Cal E_\alpha(t)$ 
satisfying 
$$
\align
q_tv_k(t) &= v_k(t), \tag{4.4.1}\\
S_k(t)^* &= v_k(t)^*\restriction p_0H, \tag{4.4.2}
\endalign
$$
for every $t>0$.  (4.4.3) implies that $S_k(t) = p_0v_k(t)$. 

We claim that 
$$
\<S_1(t), S_2(t)\>_{\Cal E_P(t)} = 
\<v_1(t), v_2(t)\>_{\Cal E_\alpha(t)}.  \tag{4.5}
$$
To see this we appeal to Proposition 1.7, which expresses 
the inner product of $\Cal E_P(t)$ in terms of the {\it minimal}
Stinespring dilation of the completely positive map $P_t$.  
We obtain such a dilation
$$
P_t(x) = V^*\pi_t(x)V, \qquad x\in \Cal B(p_0H)
$$
as follows.  

For every $x\in \Cal B(p_0H)$ let $\pi_t(x)$
be the restriction of $\alpha_t(xp_0)$ to the invariant 
subspace $K = [\alpha_t(p_0\Cal B(H)p_0)p_0H]$, and let 
$V$ be the inclusion map of $p_0H$ into $K$.  Then since
$P_t$ is the compression of $\alpha_t$ to $\Cal B(p_0H)$ 
we see that 
$$
P_t(x) = V^*\pi_t(x)V, \qquad x\in \Cal B(p_0H), 
$$
and the latter is obviously a minimal Stinespring 
representation for $P_t$.  Letting $q_t$ be the projection 
on $[\alpha_t(\Cal B(H))p_0H]$, we 
claim first that 
$$
K = \alpha_t(p_0)q_tH.  \tag{4.6}
$$
Indeed, the two projections $\alpha_t(p_0)$ and $q_t$ must
commute because $q_t$ belongs to the commutant of 
$\alpha_t(\Cal B(H))$, and 
$$
K = [\alpha_t(p_0\Cal B(H)p_0)p_0H] =
[\alpha_t(p_0)\alpha_t(\Cal B(H))p_0H] = \alpha_t(p_0)q_tH.  
$$

For $k=1,2$ we claim that the operator 
$$
X_k = v_k(t)\restriction_{p_0H}
$$
maps $p_0H$ into $K$ and satisfies 
$$
X_kx = \pi_t(x) X_k, \qquad x\in \Cal B(p_0H).  
$$
For that, note that since $v_k(t)$ belongs to 
$\Cal E_\alpha(t)$ and satisfies (4.4.1) we have 
$$
X_kp_0 = v_k(t)p_0 = q_tv_k(t)p_0 = 
q_t\alpha_t(p_0)v_kp_0 = q_t\alpha_t(p_0)X_kp_0, 
$$
and hence (4.5) implies that $X_kp_0H \subseteq K$.  Similarly,
for any operator $x$ in $\Cal B(p_0H)$ we have 
$X_kx = X_kxp_0 = \alpha_t(xp_0)X_k = \pi_t(x)X_k$.  

Finally, because of (4.4.2) we find that 
$$
S_k(t) = V^*X_k, \qquad k=1,2.   
$$
According to Proposition 1.7, the inner product 
$\<S_1(t), S_2(t)\>$ is defined by 
$$
\<S_1(t), S_2(t)\>\bold 1_{p_0H} = X_2^*X_1.  \tag{4.7}
$$
We compute the right side of (4.7).  Since 
$v_k(t)\in \Cal E_\alpha(t)$ it follows that 
$$
v_2(t)^*v_1(t) = \<v_1(t), v_2(t)\>_{\Cal E_\alpha(t)}\bold 1_H, 
$$
and thus for $\xi$, $\eta\in p_0H$, 
$$
\<X_1\xi, X_2\eta\> = \<v_1(t)\xi, v_2(t)\eta\> =
\<v_1(t), v_2(t)\>_{\Cal E_\alpha(t)}\<\xi, \eta\>.  
$$
It follows that 
$$
X_2^*X_1 = 
\<v_1(t), v_2(t)\>_{\Cal E_\alpha(t)}\bold 1_{p_0H},
$$
and (4.5) follows.  

Finally, letting $\Cal P = \{0=t_0 < t_1 < \dots < t_n=t\}$ 
be a finite partition of $[0,t]$ we find that 
$$
\align
\prod_{k=1}^n&\<S_1(t_k-t_{k-1}), S_2(t_k-t_{k-1})\> = 
\prod_{k=1}^n\<v_1(t_k-t_{k-1}), v_2(t_k-t_{k-1})\> = \\
&\<v_1(t_1-t_0)\dots v_1(t_n-t_{n-1}), 
v_2(t_1-t_0)\dots v_2(t_n-t_{n-1})\>_{\Cal E_\alpha(t)}.  
\endalign
$$
Utilizing (3.19), the last term on the right of the above 
formula is 
$$
\<q_{\Cal P, t}T_1(t), q_{\Cal P,t}T_2(t)\>_{\Cal E_\alpha(t)} =
\<Q_{\Cal P,t}T_1(t), T_2(t)\>_{\Cal E_\alpha(t)}, 
$$
and Theorem 4.3 follows \qed
\enddemo

\proclaim{Corollary 4.8} Let $P$ be a unital CP semigroup
with minimal dilation $\alpha$, and let 
$\theta: \Cal U_\alpha\to\Cal U_P$ be the bijection 
of Theorem 3.6.  Then for any two units $u_1$, $u_2$ of 
$\alpha$ we have 
$$
c_P(\theta(u_1),\theta(u_2)) = c_\alpha(u_1,u_2).  
$$
\endproclaim

\demo{proof}Let $S_i=\theta(u_i)\in\Cal U_P$, $i=1,2$.  
It is enough to show that 
$$
e^{tc_P(S_1,S_2)} = e^{tc_\alpha(u_1,u_2)}
$$
for every $t>0$.  Now Theorem 4.3 implies that 
$$
e^{tc_P(S_1,S_2)} = 
\lim_\Cal P\prod_{k=1}^n\<S_1(t_k-t_{k-1}),S_2(t_k-t_{k-1})\>
=\lim_\Cal P\<Q_{\Cal P,t}u_1(t),u_2(t)\>_{\Cal E_\alpha(t)}.  
$$
On the other hand, 
Proposition 4.2 implies that the net of projections 
$Q_{\Cal P,t}\in \Cal B(\Cal E_\alpha(t))$ increases with 
$\Cal P$ to the identity operator of $\Cal B(\Cal E_\alpha(t))$.  
Hence 
$$
\lim_\Cal P \<Q_{\Cal P, t}u_1(t), u_2(t)\>_{\Cal E_\alpha(t)} = 
\<u_1(t), u_2(t)\>_{\Cal E_\alpha(t)}.  
$$
By definition of the covariance function of $\alpha$ \cite{1}
we have 
$$
\<u_1(t), u_2(t)\>_{\Cal E_\alpha(t)} = e^{tc_\alpha(u_1, u_2)}, 
$$
as required \qed
\enddemo

With Corollary 4.8 in hand, the remarks at the beginning of this 
section imply the following,

\proclaim{Theorem 4.9}
Let $P$ be a CP semigroup and let $\alpha$ be its minimal 
dilation to an \esg .  Then 
$$
d_*(P) = d_*(\alpha).  
$$
\endproclaim

\remark{Remark 4.9}
If we are given two CP semigroups $P$ and $Q$ acting respectively
on $\Cal B(H)$ and $\Cal B(K)$, then there is a natural CP semigroup
$P\otimes Q$ acting on $\Cal B(H\otimes K)$.  For each $t\geq 0$, 
$(P\otimes Q)_t$ is defined uniquely by its action on elementary 
tensors via
$$
(P\otimes Q)_t: x\otimes y \mapsto P_t(x)\otimes Q_t(y), 
\qquad x\in \Cal B(H), y\in \Cal B(K).  
$$
Now suppose that $P$ and $Q$ are unital CP semigroups.  
Using the minimality criteria developed in \cite{3}, it 
is quite easy to see that if $\alpha$ and $\beta$ are 
respectively minimal dilations of $P$, $Q$ to \esg s acting
on $\Cal B(\tilde H)$, $\Cal B(\tilde K)$ where 
$\tilde H \supseteq H$ and $\tilde K\supseteq K$, then 
$\alpha\otimes \beta$ is a minimal dilation of the tensor
product $P\otimes Q$ to an \esg\ acting on 
$\Cal B(\tilde H\otimes \tilde K)$.   
\endremark

Thus, from Theorem 4.9 together with a) Bhat's theorem \cite{6,7}
on the existence of \esg\ dilations of CP semigroups and b)
the addition formula for the index of \esg s \cite{2}, 
we deduce
\proclaim{Corollary 4.10}
If $P$ and $Q$ are unital CP semigroups then
$$
d_*(P\otimes Q) = d_*(P) + d_*(Q).  
$$
\endproclaim


\end